\documentclass[11pt, reqno]{amsart}
\usepackage{amssymb}

\newtheorem{theorem}{Theorem}[section]

\theoremstyle{definition}

\newtheorem{definition}[theorem]{Definition}

\newtheorem{conjecture}{Conjecture}
\newtheorem{example}[theorem]{Example}

\newtheorem{question}{Question}

\def\eul{\rm{e}}
\newcommand{\divides}{|}

\def\N{\mathbb{N}}

\def\beq {\begin{equation}}
\def\endq {\end{equation}}

\renewcommand{\le}{\leqslant}
\renewcommand{\ge}{\geqslant}

\def\beq {\begin{equation}}
\def\endeq {\end{equation}}

\renewcommand{\epsilon}{\varepsilon}

\begin{document}
\title{The Repulsion Motif in Diophantine Equations}

\author{Graham Everest and Thomas Ward}
\address{School of Mathematics, University of East Anglia,
Norwich NR4 7TJ, UK} \dedicatory{To the staff of Mulbarton Ward
and the Weybourne Day Unit} \email{g.everest@uea.ac.uk}
\email{t.ward@uea.ac.uk}

\begin{abstract}
Problems related to the existence of integral and rational
points on cubic curves date back at least to Diophantus. A
significant step in the modern theory of these equations was
made by Siegel, who proved that a non-singular plane cubic
equation has only finitely many integral solutions. Examples
show that simple equations can have inordinately large integral
solutions in comparison to the size of their coefficients. A
conjecture of Hall attempts to ameliorate this by bounding the
size of integral solutions simply in terms of the coefficients
of the defining equation. It turns out that a similar
phenomenon seems, conjecturally, to be at work for solutions
which are close to being integral in another sense. We describe
this conjecture as an illustration of an underlying motif -- repulsion -- in
the theory of Diophantine equations.
\end{abstract}

\maketitle


\section{Challenging Questions}

In~$1657$, Pierre de Fermat challenged the English mathematicians Sir Kenelm Digby
and John Wallis to
find all the integer solutions to the equation
\begin{equation}\label{Fasked}
y^2+2=x^3.
\end{equation}
We can make an educated guess about his motivation.

Diophantus~\cite[Book~VI, Prob.~19]{dioph} asked
for a right triangle, the sum of whose area~$y$
and hypotenuse~$h$ is a square, and whose perimeter is
a cube. Taking the legs to be~$2$ and~$y$, and~$h+y$ to
be~$25$, he noted
that the square~$25$, when added to~$2$, gives the
cube~$27$. The two equations~$h^2=y^2+4$
and~$y+h=25$ give~$y=\frac{621}{50}$ as a solution to the
problem. This is pleasing, but we are more interested in the intermediate
observation,
namely a solution in integers to the Diophantine equation~(\ref{Fasked}).

Bachet, in his famous~$1621$ Latin translation of the
{\it Arithmetica} of Diophantus, noted that from the one
solution~$5^2+2=3^3$ other rational numbers~$r$ may be found
with the property that~$r^2+2$ is a cube. Fermat
acquired a copy of this work, and by~$1636$ had
studied it carefully and made significant
advances. He also recorded several of his
most influential marginal notes
in this book, including the comment
concerning `Fermat's last theorem'.

In~$1656$ the English mathematician Wallis, at the
time chief cryptographer to Parliament, published his
{\it Arithmetica infinitorum}. Digby brought
this to the
attention of Fermat, prompting Fermat to
begin a correspondence with `Wallis
and other English mathematicians' concerning some of his
number problems. This correspondence all passed
through Digby's hands, and thus in~$1657$ Fermat
challenged Digby and Wallis:
\begin{quote}
{\it Peut on trouver en nombres entiers un carr{\'e} autre
que~$25$, qui, augment{\'e} de~$2$, fasse un cube? A la
premi{\`e}re vue cela par{\^a}it d'une recherche difficile, en
fractions une infinit{\'e} de nombres se d{\'e}duisent de la
m{\'e}thode de Bachet; mais la doctrine des nombres entiers,
qui est assur{\'e}ment tr{\`e}s-belle et tr{\`e}s-subtile, n'a
{\'e}t{\'e} cultiv{\'e}e ni par Bachet, ni par aucun autre dont
les {\'e}crits venus jusqu'{\`a} moi.}
\end{quote}
That is,
to find all the integer solutions to the equation~(\ref{Fasked}).

The solution~$(x,y)=(3,5)$ noted by Diophantus is one, and
Digby and Wallis might have agreed it is the only one, since
negative numbers had yet to enjoy a fully equal status with
their positive siblings. Presumably the challenge was to show
that, sign apart, there are no others, as Fermat claimed a
proof that this solution is unique. It is not clear from this
distance in time if any of the three protagonists ever did have
a complete proof of this.

\subsection{What if~$2$ had been~$15$?}
What might have been the reaction if Fermat
had instead challenged the English mathematicians
with the equation
$$
y^2=x^3+15?
$$
The solution~$(1, 4)$ is easy to spot.
It is less easy to notice that~$(109, 1138)$ is also
a solution -- and not at all obvious that, issues
of sign aside, there are no others.

Even more challenging, what about the equation
$$
y^2=x^3+24?
$$
The solution~$(1, 5)$ is once again easy to spot, as is the
solution $(-2, 4)$, and it is not too difficult to find a third
solution~$(10, 32)$. However, life is surely too short to find
the solution~$(8158, 736844)$ without mechanical help. Once
again, up to sign, this is the full list of solutions.

As these examples show, finding the solutions to the
Diophantine equation~$y^2=x^3+d$ depends in a rather
unpredictable way on the constant~$d$. Among the many results
on this problem, some stand out. Euler showed that there are no
positive rational solutions apart from~$(2,3)$ when~$d=1$ using
the method of descent, and a long series of other special cases
were solved by many other mathematicians (see
Dickson~\cite[Chap.~XX]{MR0245500} for the details). A snapshot
of the state of knowledge on this question prior to Siegel's
theorem may be found in the work of
Mordell~\cite{jfm.44.0230.03}, where many but not all values
of~$d$ for which there are no integral solutions are found.

Nowadays the numerical facts above -- not only the stated
solutions, but the much deeper claim that there are no other
integral solutions -- can be checked easily using one of the
many sophisticated computational packages available, such as
{\sc Magma}~\cite{MR1484478}. This is greatly to be celebrated,
although the facility with which these calculations can now be done
risks obscuring the remarkable achievements made in Diophantine
analysis over the last forty years that have made this
possible. It is not our intention to survey these achievements
(a detailed overview may be found in the lovely monograph
of Hindry and Silverman~\cite{MR1745599}), although some
highlights on the theoretical side will appear naturally.

\subsection{Satisfying Fermat} Out of a patriotic desire to satisfy Fermat,
we start with a (now elementary) proof that~$(3,\pm 5)$ are
indeed the only integral solutions to~$y^2=x^3-2$. The language
of (and one result from) modern algebra makes this
straightforward. Plainly~$x$ must be odd, for otherwise~$y^2
\equiv 2\pmod{4}$, which is impossible. The ring~$R=\mathbb
Z[\sqrt{-2}]$ is a unique factorization domain (indeed, is a
Euclidean domain -- see Hardy and
Wright~\cite[\S14.7]{MR568909} for an account) -- whose only
units are~$\pm 1$. Factorizing there gives
\[
(y+\sqrt{-2})(y-\sqrt{-2})=x^3.
\]
In the ring~$R$,~$\gcd(y+\sqrt{-2},y-\sqrt{-2})\divides 2\sqrt{-2}.$
Since $x$ is odd, the greatest common divisor
is a unit, and so must be~$\pm 1$.
It follows that each factor is a unit multiple of a cube in~$R$.
Write (absorbing~$-1=(-1)^3$ if necessary)
\[
y+\sqrt{-2} = (a+b\sqrt{-2})^3
\]
with~$a,b\in\mathbb Z$.
Comparing coefficients of~$\sqrt{-2}$ gives
\[
1=3a^2b-2b^3=b(3a^2-2b^2)
\]
which forces~$b=3a^2-2b^2=\pm 1$.
Solving gives~$b=1$ and~$a=\pm1$, and
the two choices give~$y=\mp 5$ and~$x=3$.

\section{Siegel's Theorem}
The following is a special case of a wonderful and deeply
influential result, proved by Siegel. Siegel's
result was more general, but for the case at hand may be stated
as follows.

\begin{theorem}[Siegel~\cite{siegel1929}, 1929]\label{siegel}
Given an integer~$d\neq 0$, the equation
\begin{equation}\label{mordellequation}
y^2=x^3+d
\end{equation}
has only finitely many integral solutions.
\end{theorem}

This specific result was known earlier.
Mordell~\cite{mordell1920} reported to the London Mathematical
Society in~$1918$ that an earlier result of
his~\cite{mordell1914}, in conjunction with Thue's
work~\cite{thue1909}, would show Theorem~\ref{siegel}. 
Once again there is an important letter, sent by Siegel to
Mordell in~$1925$ (an extract appears in~\cite{siegel1926}),
outlining Siegel's ideas which eventually gave the general
result~\cite{siegel1929} that any non-singular cubic curve has
only finitely many integer points. In honour of Mordell's contribution
to this subject, the equation (\ref{mordellequation}) in
Theorem~\ref{siegel} is often known as~\emph{Mordell's
equation}.

Even in its simplest form, Siegel's proof was
recognized as \emph{non-effective}: a term
whose meaning will be discussed shortly. Siegel gave a second 
proof of finiteness using unique factorization, in a
manner very similar to the proof in the case~$d=-2$ above, by
working in a suitably large ring and then reducing to a number
of simpler equations called~$S$-unit equations. Later this led
to an \emph{effective} proof, following Alan Baker's seminal
work on transcendence theory. The terms effective and
non-effective, although ubiquitous in number theory, are never
precisely defined.

\begin{quote}
ef¸fec¸tive {\it
(adj.):
adequate to accomplish a purpose; producing the intended or expected result;}
1350--1400; ME fr. L. {\it effectivus} $=$ practical.
\end{quote}
In a nutshell, still lacking a mathematical
definition, one might
expect that an
effective proof is one which produces an algorithm to implement the
conclusion of a theorem.
For example,
in our context, an effective proof might consist of a
bound on the size of the largest solution.
This would allow -- \emph{in principle} -- all solutions
to be found, by simply
checking integers below that bound to see if they
satisfy the equation.
This sounds desirable and practical, but there are examples
where the gap between
what theory provides and what is practical remains large.

\begin{theorem}[Baker~\cite{MR0231783}, 1968]\label{baker}
Any integral solution of {\rm(\ref{mordellequation})} satisfies
\[
|x|<\eul^{10^{10}|d|^{10^4}}.
\]
\end{theorem}

How large is the number on the right hand side? It is said that
the number of electrons in the visible universe does not
exceed~$10^{80}$, while those on the earth amount to a little
over~$10^{60}$. Each of these quantities is minuscule compared
to the bound in Theorem~\ref{baker}, even for modest values of~$d$.
It is tempting to imagine
that with the speed of modern computers, it would nonetheless
be practical to check values of~$x$ below that bound, giving a
complete solution for any given~$d$
of reasonable size. Before attempting this, it
is wise to think through some realistic estimates for the type
of computer one might have at home. When~$d=-2$, you will
probably find it takes under a minute to check if~$x$
satisfies~\eqref{mordellequation} for~$\vert x\vert<10^6$, over
an hour to reach~$10^9$. There are roughly~$10^9$ personal
computers in the world (and this figure is likely to double
by~$2015$). If you owned them all and used them in your
checking, it would take at least a billion billion billion
billion billion billion years to check the equation for values
of~$x$ with~$|x|<10^{80}$, and this is well short of the bound
in Theorem~\ref{baker}. Even checking Euler's proof that~$x=2$
gives the only positive integral solution to~$y^2=x^3+1$ is not
feasible using Baker's bound.

So, despite the etymology, \emph{effective} as used in number
theory (which would certainly include the statement in
Theorem~\ref{baker}) does not always mean practical.

Of course Theorem~\ref{baker} was never intended to be a
practical tool, and was never claimed to be such. It merely
says that with the methods then available, this is the best
bound obtainable, and sets a challenge for future generations.
Dramatic improvements followed fairly quickly. For example,
Stark~\cite{MR0340175} shows that
\[
\max\{
\vert x\vert,\vert y\vert\}
\leq
\eul^{C\vert d\vert^{1+\epsilon}},
\]
where~$C=C(\epsilon)$ is an effectively computable constant
depending only on~$\epsilon$. This is a considerable
improvement upon Baker's bound, however, the size
of~$C(\epsilon)$ is necessarily large, and checking small
values is not computationally feasible, even for small values
of~$|d|$.

Despite the power of these results and the direction of
research initiated by the work of Baker and Stark, the modern
computational facility for solving Mordell's equation came
about \emph{not} by further reducing the size of the bounds
above. The method is actually less direct (see the work of
Gebel, Peth{\"o}, and Zimmer~\cite{MR1464020} for complete
details about the toolkit now used). In a later
paper~\cite{MR1602064}, the authors show how this method
resolves the equation for~$|d|\le10^4$, and for almost all~$d$
with~$|d|\le10^5$.

Today, after major theoretical and computational improvements,
computer packages will now find all the integral solutions of
equations~$y^2=x^3+d$ provided~$d$ lies within
\emph{reasonable} bounds. To give an idea of what counts as
\emph{reasonable} it is worthwhile doing some experiments
yourself.

\section{Hall's Conjecture} To reiterate, solving Mordell's equation in practice does not rely
upon obtaining a very strong upper bound for the size of the
largest solution, then checking smaller solutions on a
computer. In fact an important, simple, natural, question
remains unsolved to this day.

\begin{question}
What is the best theoretical bound (in terms of~$d$) for the
size of the largest integer~$x$ solving~$y^2=x^3+d$?
\end{question}

The truth is that the best known bound is far from what seems
likely to be true. What seems likely to be true is the subject
of the following conjecture made by Hall~\cite{MR0323705},
which we will now discuss. Hall's conjecture has been subject to
extensive numerical checking, but a proof seems to require
dramatically more powerful methods than those currently
available.

\begin{conjecture}[Hall\footnote{Originally, Hall \cite{MR0323705} conjectured
the same bound but with $\epsilon =0$. This is no longer
thought to be likely.}]
Given~$\epsilon>0$, there is a constant~$C=C(\epsilon)$ such
that, for any non-zero~$d\in\mathbb Z$, any integral solution
of~$y^2=x^3+d$ satisfies
\[
\log\vert x\vert<(2+\epsilon)\log\vert d\vert+C.
\]
\end{conjecture}

The audacious nature of the conjecture is not immediately
apparent. Our second and third examples ($d=15$ and~$d=24$)
show that simple equations with small coefficients can have
enormously large integral solutions. What Hall's conjecture
suggests is that, when properly calibrated, the phenomenon of
large integral solutions of an equation with small coefficients
is not beyond constraint.

The conjecture of Hall follows from the infamous~\emph{$ABC$
conjecture}, formulated by Masser and
Oesterl{\'e}~\cite{MR992208} in~$1985$, about the relative
sizes of a zero sum of three integers. Write
\[
r(N)=\prod_{p|N}p
\]
for the \emph{radical} of an integer~$N$. The~$ABC$ conjecture
says that for any~$\epsilon>0$ there is a
constant~$C=C(\epsilon)$ such that, whenever
\[
A+B+C=0
\]
in non-zero coprime integers~$A,B,C$, we have
\[
\max\{|A|,|B|,|C|\}\le Cr(ABC)^{1+\epsilon}.
\]

To see how this relates to the Mordell equation, assume
that~$|x|^3$ is the largest of the three terms
in~\eqref{mordellequation}. Then, roughly speaking, the~$ABC$
conjecture predicts that
\[
|x^3|\le C |xyd|^{1+\epsilon},
\]
and the Hall bound follows by taking logs and noting
that~$|y|$ is approximately~$|x|^{3/2}$.

Hall's conjecture has been extensively tested.
Table~\ref{tableonSHconjecture} shows values of integers~$x$ and~$d$
with an integral~$y$ satisfying~$y^2=x^3+d$, having~$\log x $ large in comparison
with~$2\log\vert d\vert$. It is taken from Elkies'
website~\cite{elkieshall} (see also the
paper~\cite{MR1602064}). The table is surprising in two
opposite senses: Firstly, it gives more examples of
inordinately large solutions of simple Diophantine equations,
and secondly, it shows that they nonetheless fall within sight
of a reasonable constraint upon how large they could be
when viewed on a
logarithmic scale.

\begin{center}
\begin{table}[ht]
\caption{\label{tableonSHconjecture}Large values of~$\log
x/2\log |d|$ in Hall's conjecture.}
\begin{tabular}{|l|l|l|c|}
\hline &&&\\[-10pt]
$d$ & $x$ & $\log x$& $\log x/2\log |d|$\\
\hline
\small{-1641843}& \small{5853886516781223}& \small{36.305}
& \small{1.268}
\\
\hline
\small{-30032270}& \small{38115991067861271}&\small{38.179}
&\small{1.108}
\\
\hline
\small{1090}&\small{28187351}&\small{17.154}
&\small{1.226}
\\
\hline
\small{193234265}&\small{810574762403977064}& \small{41.236}
&\small{1.080}
\\
\hline
\small{17}&\small{5234}&\small{8.562}
&\small{1.511}
\\
\hline
\small{225}&\small{720114}&\small{13.487}
&\small{1.245}
\\
\hline
\small{24}&\small{8158}&\small{9.006}
&\small{1.417}
\\
\hline
\small{-307}&\small{939787}&\small{13.753}
&\small{1.200}
\\
\hline
\small{-207}&\small{367806}&\small{12.815}
&\small{1.201}
\\
\hline
\small{28024}&\small{3790689201}&\small{22.055}
&\small{1.076}
\\
\hline
\end{tabular}
\end{table}
\end{center}

Does Table~\ref{tableonSHconjecture} convince? Taken together with the implication
of the~$ABC$ conjecture, the answer is probably yes, in the
sense that this is evidence for a sensible conjecture. We have
presented it because, later on, two more tables will appear and
a direct comparison will be invited.

\section{Repellent Powers}\label{repellentpowers}

It might not be apparent so far, but the point at infinity is
enormously important in understanding solutions of Mordell's
equation. For example, the set of rational points on the curve
forms a group under a natural geometric form of addition. The
point at infinity is the identity element for this group.

From our viewpoint, Siegel's Theorem may be interpreted to say
that the point at infinity \emph{repels} the integer points on
the curve~$y^2=x^3+d$. In other words, there is a (punctured)
neighbourhood of the point at infinity free of integer points.
Empirically, we might even say we observe
integer points repelling each other. 
For a sophisticated instance of this repulsion property being used
to understand integral points on elliptic curves, see the work
of Helfgott and Venkatesh~\cite{MR2220098}. In their work, integer
points repel each other inside the Mordell-Weil lattice. The methods
of~\cite{MR2220098} will yield explicit constants that might
well quantify practically the rate at which integral points
repel each other coordinate-wise. This repulsion between integral points will
be something of a mantra throughout this paper, and will inform
the latter part significantly. For now though, consider the
idea of points with fixed arithmetic properties repelling each
other as a kind of paradigm for understanding other results in
Diophantine equations.

The results stated so far may be seen as an instance of a
general tendency for distinct \emph{integral powers} to repel
each other. Thus, for example, Baker's result may be phrased as
follows. If~$x$ and~$y$ are positive integers with~$y^2\neq
x^3$ then there is a constant~$C=C(x)$ with
\[
\vert y^2-x^3\vert>C(\log x)^{10^{-4}};
\]
Stark's result says that for
any~$\kappa<1$ there is a constant~$C=C(\kappa)$
with
\[
\vert y^2-x^3\vert>C(\log x)^{\kappa},
\]
and Hall's conjecture says that for
any~$\epsilon<1/2$ there is a constant~$C=C(\epsilon)$ with
\[
\vert y^2-x^3\vert>Cx^{\frac12-\epsilon}.
\]

These are all statements about (subsequences of) the sequence
\[
a=(a_n)=(1,4,8,9,16,25,27,32,36,49,64,81,100,121,125,\dots)
\]
of perfect powers: numbers of the form~$n^m$ with~$n,m\in\N$
and~$m\ge2$. A gap of one is seen early on, and the conjecture
of Catalan~\cite{catalan} is that~$8$ and~$9$ are the only
consecutive pair in the sequence~$a$. The proof of this result
followed a path that once again illustrates the slightly
ambiguous way in which the word {\it effective} is used.
Tijdeman~\cite{MR0404136} used sharpened versions of Baker's
theorem to find a number~$T$ with the property that any
positive integral solution~$(x,y,m,n)$ to the Diophantine
equation~$x^m-y^n=1$ must have~$\max\{x,y,m,n\}\le T$. This
meant that a proof of Catalan's conjecture was reduced to a
finite list of possibilities to check -- surely the most
effective of effective statements. Unfortunately a by-product
of the transcendence methods used was that the number~$T$ was
enormous, leaving a finite but hopelessly impractical
calculation to be done. The conjecture was finally proved by
Mih\u{a}ilescu~\cite{MR2076124} using a mixture of methods, and
eventually a beautiful, purely algebraic, proof was given by
Bilu~\cite{MR2152212}.

Pillai studied many properties of the sequence of perfect
powers, and in a paper~\cite{MR0013386} of~$1945$ wrote
\begin{quote}
\it I take this opportunity to put in print a conjecture which
I gave during the conference of the Indian Mathematical Society
held at Aligarh. Arrange all the powers of integers like
squares, cubes, etc. in increasing order [...]. Let~$a_n$ be
the~$n$th member of this series [...]. Then
\[
\liminf_{n\to\infty}(a_n-a_{n-1})=\infty.
\]
\end{quote}
The conjecture of Pillai is exactly equivalent to the
conjecture that for any~$k\ge1$ the Diophantine
equation~$x^m-y^n=k$ has only finitely many solutions. This
remarkable problem remains open, and is the subject of a recent
survey by Waldschmidt~\cite{pillai}.

\subsection{The Gap Principle}

The phenomena we are describing in this paper will be
recognizable to workers in Diophantine equations although
possibly under
another name -- the gap principle. Again, this is not
formulated precisely anywhere that we can find. Roughly stated,
it says that where Diophantine phenomena occur (say, as
rational solutions to a Diophantine equation or to an
inequality) subject to some reasonable constraint, they will
respond by exhibiting measurable gaps. This data is then fed
back into the technicalities of the argument. Strictly
speaking, there is no single gap principle, it is more of a
style of argument. We invite readers to explore the current
literature on Diophantine equations to see where this term --
or this phenomenon -- occurs.

We cite two examples, one old and one recent and very germane: Ingram~\cite{MR2468477}
showed that large gaps
occur between integral multiples of points on Mordell curves in the following sense. If~$n_1<n_2$
and~$n_1P,n_2P$ (in the group theory sense) are both integral points, then~$a^{n_1^2}<n_2$
for some explicit constant~$a>1$. This is a strong repulsion
property for integer points along a sequence~$(nP)$ and, just as with the
results in \cite{MR2220098}, might well translate back into a good
bound for the corresponding coordinates. An instance of the kind of
applications Ingram obtains is that on Mordell curves of the
form~(\ref{mordellequation}), for all large enough~$d$ with~$d$
sixth power-free (that is, not divisible by~$a^6$ for any integer~$a\ge2$),
there is at most one~$n\ge 3$ such that~$nP$
is integral (see~\cite[Proposition~15]{MR2468477}).

Perhaps the earliest -- and certainly a very influential --
manifestation of a gap principle occurs in an old paper of
Mumford~\cite{MR0186624}. His result was about rational points
lying on more complicated plane curves. Complicated here means
not just higher degree but higher \emph{genus}, a geometric
measure of complexity -- for example, Mumford's results apply
to equations~$y^2=f(x)$ where~$f(x)\in \mathbb Q[x]$ has degree
at least~$5$, and~$f(0)\neq 0$. He showed that the rational
solutions exhibit naturally occurring, expanding, gaps (when
viewed in the most na\"ive sense, so that the numerators and
denominators grow large very quickly). Although in one sense
his result was superseded by Faltings' general proof that such
curves contain only finitely many rational points, packages such as {\sc Magma}~\cite{MR1484478}
will now enumerate rational points on higher genus curves with
some ease. The repelling nature of them brings Mumford's remark
vividly to life.

\begin{example}[Taken from \cite{BT}] The genus 2 curve $y^2=x^6+1025$ has the following
rational points, together with sign changes:
$$(2, 33), \left(\frac{5}{2}, \frac{285}{2^3}\right), (8, 513), \left(\frac{1}{4}, \frac{2049}{4^3}\right),
\left(\frac{20}{91}, \frac{24126045}{91^3}\right).
$$
\end{example}

\section{Generalizing Siegel's Theorem and Hall's Conjecture}

We will now describe a recent attempt to generalize both Siegel's
Theorem and Hall's conjecture in one go. We need to begin by talking
a little about rational solutions of our equations. Although it is not obvious, 
there are infinitely many
rational solutions to each of our three starting equations. In
each case, there is a way to produce them all, starting with a
finite set of rational points. This fact, casually stated,
obscures the enormous theoretical and practical knowledge about
rational solutions (as opposed to integral solutions). For
example, the rational points on
\begin{equation}\label{15}
y^2=x^3+15.
\end{equation}
are all generated from the two points~$(1, 4)$
and~$(\frac{1}{4},\frac{31}{8})$, using the chord and tangent
method of constructing new points. The essential observation is
that the line joining any two rational points on the curve (or
the tangent to a single rational point)
intersects the curve again at a third rational point. For
example, the line joining our two initial points meets the
curve again at~$(-\frac{11}{9},\frac{98}{27}).$

This operation, together with reflection in the~$x$-axis,
allows all rational points on the curve to be found --- but
this is not an easy result to prove. To help orient the reader,
notice that this is an instance of \emph{Mordell's
Theorem}~\cite{mordell1922}, which says that the set of
rational solutions of such an equation is always generated from
a finite set of rational solutions in the same geometric way.

Indeed, the operation of joining two points with a line,
finding the third point of intersection, and then reflecting in
the~$x$-axis (extended by continuity to allow the original
points to be identical) is a binary operation giving the set of
rational points on the curve the structure of an abelian group
(once the point at infinity is added),
and Mordell's theorem states that this group is finitely
generated. It must therefore have the form~$\mathbb Z^r\times
F$ for some finite group~$F$, and~$r$ is called the \emph{rank}
of the curve. The quantity~$r$, like~$d$
in~\eqref{mordellequation} or~$\Delta_E$ in
Section~\ref{sectioncomputationalevidence}, is an important
measure of the size or the complexity of the curve.

To fix notation both now and for the sequel, note that
if~$P=(x,y)$ is a rational solution of~$y^2=x^3+d$ and~$d$
is integral, then~$x^3$
and~$y^2$ have the same denominator. Thus the denominator
of~$y^2$ must be the square of a cube, and that of~$x^3$ the
cube of a square, so we may write
\[
P=\left(\frac{A_P}{B_P^2}, \frac{C_P}{B_P^3}\right),
\]
with~$A_P, B_P, C_P\in \mathbb Z$ in lowest terms. 

Siegel's
Theorem makes a statement about the rational solutions~$P$ of
the equation~$B_P=1$, which we may interpret as a statement
about the solutions~$(x,y)$ under the condition that the
denominator of~$x$ (or of~$y$) is divisible by no primes.
\begin{question}
What can be said about the points~$P$ such that~$B_P$ is
divisible by one prime?
\end{question}
In other words, what can be said about the set of rational
solutions~$P$ with the property that~$B_P$ is a prime power?
There is a name for points of this form. If the point~$P$
has~$B_P$ composed of primes from a set~$S$ then~$P$ called
an~\emph{$S$-integral point}. The Siegel-Mahler Theorem
predicts that for a {\it fixed} finite set of primes $S$, there are only finitely
many $S$-integral points. What we are doing is slightly different,
in fixing the size of $S$
(typically $|S|=1$), but not its contents.
Reynolds~\cite{reynolds} proved (under a general form of
the~$ABC$ conjecture) that only finitely many pure powers will
occur among the~$B_P$, so the first interesting case is when~$B_P$ is a prime.

Sometimes, it is provable that only finitely many points~$P$
have~$B_P$ equal to a prime. Indeed, this is true of our
starting equation
\begin{equation}\label{-2}
y^2=x^3-2.
\end{equation}
To see this, note that on the curve~\eqref{-2}, all the rational points can be
generated by the single point~$(3,5)$. Since this point is the
image of a rational point under a 3-isogeny, \cite[Theorem 1.3]{MR2045409}
applies to demonstrate the finiteness claim. Inspection suggests that, in truth,
none of the rational points yield prime values~$B_P$ but, in a familiar refrain,
the amount of checking of small cases remains, as yet, unfeasible.

On the other hand,
searching yields a large number of rational points~$P$
on equation~\eqref{15} with~$B_P$ equal to a prime. The examples below,
and all that follow, were obtained using searches in
Pari-GP~\cite{parigp}. What prevents us exhibiting many more examples
is lack of space, not data. Readers interested in following this up should
consult~\cite{MR2041076}.

\begin{example} The following~$x$ coordinates of rational
points~$P$ on the curve~$y^2=x^3+15$ have~$B_P$ equal to a
prime:
\begin{eqnarray*}
x(P)&=&-\textstyle\frac{11}{3^2},\\
x(P)&=&\textstyle\frac{75721}{53^2},\\
x(P)&=&\textstyle\frac{578509}{367^2},\\
x(P)&=&-\textstyle\frac{349755479}{11909^2},\\
x(P)&=&\textstyle\frac{556386829130869}{17684189^2},\mbox{ and}\\
x(P)&=&\textstyle\frac{64892429414388628056900713281}{259476976750177^2}.\\
\end{eqnarray*}
\end{example}

\section{The condition that $B_P$ is a prime}

Computational evidence, as well as a heuristic
argument~\cite{MR2045409,MR2041076}, suggest that when the set
of rational points is generated by more than one point, as in
equation~\eqref{15} for example, there will be infinitely many
points~$P$ with~$B_P$ equal to a prime\footnote{Provided a
technical assumption is met, namely, that the set of points
should not lie in the image of a rational isogeny.}. Although
the papers~\cite{MR2045409,MR2041076} contain a great deal of
numerical evidence, as well as a reasoned heuristic argument,
we must admit that the conjecture has not been proved, even for
one curve.

The question we now ask can be put roughly as follows.

\begin{question}
Where are the rational points~$P$ with~$B_P$ a prime?
\end{question}

It makes sense to set the question in a slightly wider
context, even though our main concern lies with the case
stated.

\begin{definition}
The \emph{length} of a rational point~$P$ is the number of
distinct prime divisors of~$B_P$.
\end{definition}

Thus the integral points are precisely the length zero points.
Length~$1$ points are those with~$B_P$ equal to a prime power.
Our comments will apply to rational points whose length is
bounded, but the question above asks simply for the location of
the length~$1$ points. As discussed in
Section~\ref{repellentpowers}, Siegel's Theorem may be
interpreted to say that the point at infinity repels length
zero points. Is it possible that the same might be true for
length~$1$ points? In other words, do length~$1$ points have
bounded~$x$-coordinates?

A first attack is computational. To increase the options, a
search was made on curves~$y^2=x^3+Ax+B$ with~$4A^3+27B^2\neq
0$, that is, on elliptic curves (which is the proper setting
for this question for reasons we will expand on later). This
does yield examples of length~$1$ points with inordinately
large~$x$-coordinates.

\begin{example}\label{whoops}
The curve~$y^2=x^3-7x+10$ contains a rational point~$P$ with
\[
B_P=14476032998358419473538526891666573479317742071,
\]
a prime, and with
\[
x(P)= 175567.984\dots.
\]
\end{example}
The $x$-coordinate has been expressed as a real number to emphasize
its rough size -- this is less apparent when the number is written in
rational form with a numerator and denominator.

\subsection{Generalizing Siegel's Theorem}
Does this mean that the obvious generalization of Siegel's Theorem is false?
One might think so given
Example~\ref{whoops} above, and others like it -- but that is to miss the flow
of our argument thus far.
Drawing a parallel with our earlier comments, the evidence both from examples of
Siegel's Theorem and from the constraints along the lines of
Hall's Conjecture suggests that the right parameter to measure
is~$\log x /\log |d|$ (for a suitable notion of the
parameter~$d$ adapted to more general curves). Viewed on this
logarithmic scale, we present computational evidence (see Table~\ref{ellipticcurves})
to suggest that
infinity does indeed repel length~$1$ points. In other words, we are suggesting that the
generalization of Siegel's Theorem becomes more
plausible when strengthened
along the lines of Hall's conjecture.

More even than this, an elliptic curve is fundamentally a projective object
(in particular,
the process of adding the point at infinity as the identity
element can be formalized by viewing the curve as an elliptic object).
What this means is that its true nature only becomes revealed when viewed
as lying in projective space. On that basis, there are no specially favoured
points: although the point at infinity is chosen traditionally as the identity for the
group law, the group structure makes sense with any chosen rational
point as identity, with the appropriate changes. Thus we are drawn, perhaps with some trepidation,
towards a belief that
\emph{all} rational points (or even all algebraic points)
repel the points with length below a fixed bound. In other words, around each algebraic point there is a punctured
neighbourhood free of bounded length points.

In order to examine this conjecture, as well as to frame it along earlier lines, we
introduce a suitable notion of the distance between two points.

\begin{definition}\label{logdist}
Given an algebraic point~$Q$, we define the \emph{logarithmic
distance} from~$P \neq Q$ to be
\[
h_Q(P)=-\log |x(Q)-x(P)|
\]
when~$Q$ is finite, and
\[
\log |x(P)|
\]
when~$Q$ is the point at infinity.
\end{definition}

To help swallow the looking-glass nature of
Definition~\ref{logdist}, note that points close to `infinity'
are considered large. For the sake of consistency, logarithmic
distances work the same way, so that points which are close to
each other turn out to have a large logarithmic distance -- a
general feature of such distances in Diophantine Geometry.
To see why this is natural, notice that a measure of the
size of a rational number~$\frac{a}{b}$ in lowest terms
is given by the \emph{height}~$H(\frac{a}{b})=\max\{\vert a\vert,\vert b\vert\}$.
Then a rational number very close to, but not equal to, and
specified number necessarily has very large height.

\begin{conjecture}\label{newhall} Let $k$ denote a fixed positive integer. Assume $d$
is sixth power-free. Then the set of algebraic points repels the
rational points with length below~$k$ (the `Siegel part' of the
conjecture). If~$Q$ denotes a fixed algebraic point, then there
is a constant~$C=C(k,Q)$ such that
\[
h_Q(P)\le C\log |d|
\]
for any rational point~$P$ with length below~$k$ (the `Hall
part').
\end{conjecture}

The condition on~$d$, that it is not divisible by any sixth
power, is a natural one. Without this condition, we would be
free to multiply the equation through by sixth powers of
integers to create more and more integral and
length~$1$ points in an
artificial way. The general hypothesis, of which this is the
special case for~$y^2=x^3+d$, is that an elliptic curve be in
\emph{minimal form}.

In the next section, computational evidence for
Conjecture~\ref{newhall} will be presented. Notice that some
evidence for Conjecture~\ref{newhall} already comes from the
data examined in the first part of the article. The conjecture
predicts that integral points are not simply repelled by the
point at infinity, but also by each other -- exactly what we
observe. Conversely, the conjecture too sheds light upon the
data. Looked at this way, we should not be surprised that, when
more than a couple of integral points exist, the outliers are
forced to lie a long way out.

\section{Computational
Evidence}\label{sectioncomputationalevidence}

This was obtained in \cite{MR2548983} for elliptic curves
(see~\cite{MR1144763,MR2514094} for background) of the form
\[
E: y^2+a_1xy+a_3y=x^3+a_2x^2+a_4x+a_6
\]
with~$a_1,\dots ,a_6\in \mathbb Z$. The role of the
parameter~$d$, measuring the complexity of the curve, is played
here by~$0\neq\Delta_E\in\mathbb Z[a_1,\dots ,a_6]$, the
so-called \emph{discriminant}. This is a (complicated) explicit
polynomial in the variables~$a_1,\dots,a_6$. Write~$h_E=\log
|\Delta_E|$. Each curve~$E$ is recorded as a
list~$[a_1,a_2,a_3,a_4,a_6]$, and Table~\ref{ellipticcurves}
shows what we believe is the maximal distance $\overline h$
from infinity of the length~$1$ points on that curve.

The curves are all in minimal form and all have rank~$2$, and
the search range runs over length~$1$ points of the
form~$mP+nQ$ where~$P$ and~$Q$ denote generators of the
torsion-free part of the group of rational points, and~$|m|,
|n|\le150$. The search range is necessarily constrained because
checking that points have length~$1$ requires primality testing on
some large integers. In every case, the largest
value~$\overline h$ occurs fairly early, strengthening our
belief it is truly the maximum -- see~\cite{MR2548983} for
details. Conjecture~\ref{newhall} predicts that the
ratio~$\overline{h}/h_E$ will be uniformly bounded. Readers are
invited to compare this table of values with the one for Hall's
Conjecture itself (Table~\ref{tableonSHconjecture}) given
earlier.

\begin{center}
\begin{table}[h]
\caption{\label{ellipticcurves} Infinity repelling length $1$ points}
\begin{tabular}{|l|l|l|l|l|}
\hline &&&\\[-10pt]
$E$  & $|\Delta_E|$ & $\overline h$ &$\overline h/h_E$ \\[3pt]
\hline
[0,0,1,-199,1092]&
11022011& 12.809 & 0.789 \\
\hline
[0,0,1,-27,56]&107163&11.205& 0.967\\
\hline
[0,0,0,-28,52]&
236800&
13.429&
1.085\\
\hline
[1, -1, 0, -10, 16]&
10700&
9.701&
1.045\\
\hline
[1,-1,1,-42,105]&
750592&
8.136&
 0.601\\
\hline
[0, -1, 0, -25, 61]&
154368&
16.592&
1.388\\
\hline
[1, -1, 1, -27, 75]&
816128&
 12.363& 0.908\\
\hline
[0, 0, 0, -7, 10]&
21248 &
12.075& 1.211\\
\hline
[1, -1, 0, -4, 4]&
892&
11.738& 1.727\\
\hline
[0, 0, 1, -13, 18]&
3275&
6.511& 0.804\\
\hline
[0, 1, 0, -5, 4]&
4528&
7.377& 0.876\\
\hline
[0, 1, 1, -2, 0]&
389&
9.707& 1.627\\
\hline
[1, 0, 1, -12, 14]&
2068&
9.819& 1.286\\
\hline
\end{tabular}
\end{table}
\end{center}

Table~\ref{secondellipticcurvetable} shows some curves
with~$Q=(0,0)$, and the maximal distance~$\overline h_Q$
from~$Q$ of the length~$1$ points. The remarks about search
ranges and our confidence about the true nature of the maximum
apply here as before, see also~\cite{MR2548983}. Note that
points~$P$ close to~$Q=(0, 0)$ yield a large values of~$h_Q(P)$
because these are logarithmic distances. The thrust of
Conjecture~\ref{newhall} is that these values, properly scaled,
are not inordinately large, but satisfy a reasonable constraint.
Once again, a direct comparison is invited between this table
and the one given earlier in Table~\ref{tableonSHconjecture}.

\begin{center}
\begin{table}[h]
\caption{\label{secondellipticcurvetable}Curves
with~$Q=(0,0)$.}
\begin{tabular}{|l|l|l|l|}
\hline &&&\\[-10pt]
$E$  & $|\Delta_E|$  & $\overline h_{Q} $ &$\overline h_{Q}/h_E$ \\[3pt]
\hline
[0, 0, 0, 150, 0] &216000000&  6.436 & 0.335\\
\hline
[0, 0, 0, -90, 0] & 46656000&   3.756 & 0.212\\
\hline
[0,0,0,-132,0] &147197952 &4.470&0.237\\
\hline
[0,1,0,-648,0] & 17420977152& 0.602&0.025\\
\hline
[0,0,0,34,0] &2515456& 2.107&0.143\\
\hline
[0,0,0,-136,0] & 160989184& 0.279&0.014\\
\hline
[0,1,0,-289,0] & 1546140752 &5.712&0.269\\
\hline
\end{tabular}
\end{table}
\end{center}

\section{A Theorem}

Given the speculative nature of this study, it is a little
surprising to find certain conditions where the Siegel part of
the conjecture can be proved unconditionally, and the Hall part
conditionally, at least when~$Q$ is the point at infinity. What
follows is a modest result, but it does yield examples where further
testing may be carried out, and it provides some support for
Conjecture~\ref{newhall}. Recall that rational points~$Q_1$
and~$Q_2$ are \emph{independent} in the group law on the
rational points of the curve if there is no point~$P$
with~$Q_1=aP$ and~$Q_2=bP$ for integers~$a,b$.

\begin{theorem}[Everest and
Mah{\'e}~\cite{MR2548983}]\label{gmtheorem} Consider the
equation
\[
y^2=x^3-Nx
\]
for some~$0<N\in\mathbb N$. Assume that~$Q_1$ and~$Q_2$ are
independent rational points
with~$x(Q_1)<0$ and with~$x(Q_2)$ equal to a square.
\begin{enumerate}
\item (`Siegel part') There is a bound upon $|x(P)|$ as $P$ runs
over length~$1$ points in the
    group generated by~$Q_1, Q_2$.
\item (`Hall part') Assume additionally that the~$ABC$
    conjecture holds in~$\mathbb Z$. Then there is a
    constant~$C$, independent of~$N$, for which
\[
    \log |x(P)|\le C \log N,
\]
for all length $1$ rational points~$P$.
\end{enumerate}
\end{theorem}

\begin{example}
The points~$Q_1=[-9,9], Q_2=[49/4,-217/8]$ on the curve
\[
y^2 = x^3 - 90x
\]
satisfy the hypotheses of Theorem~\ref{gmtheorem}. Computations
support the belief that infinitely many length~1 points lie in
the group generated by~$Q_1$ and~$Q_2$.
\end{example}

\begin{example}
The points~$Q_1=[-9,120], Q_2=[841,24360]$ on the curve
\[
y^2=x^3-1681x
\]
satisfy the hypotheses of Theorem~\ref{gmtheorem}. Note
that~$x(Q_2)=29^2$. Also, in this case~$Q_1$ and~$Q_2$ are
generators for the torsion-free part of the group of rational
points. As before, it seems likely that infinitely many length~1 points lie in
the group generated by~$Q_1$ and~$Q_2$.
\end{example}


\begin{thebibliography}{10}

\bibitem{MR0231783}
A.~Baker,  `The {D}iophantine equation {$y^{2}=ax^{3}+bx^{2}+cx+d$}', \emph{J.
  London Math. Soc.} \textbf{43} (1968), 1--9.
\newblock \verb|http://dx.doi.org/10.1112/jlms/s1-43.1.1|.

\bibitem{parigp}
C.~Batut, K.~Belabas, D.~Bernardi, H.~Cohen, and M.~Olivier, \emph{User's guide
  to {PARI-GP}} (Laboratoire A2X, Universit\'e Bordeaux I, France, 1998).
\newblock
  \verb|http://pari.math.u-bordeaux.fr/pub/pari/manuals/2.1.6/users.pdf|.

\bibitem{MR2152212}
Y.~F. Bilu,  `Catalan without logarithmic forms (after {B}ugeaud, {H}anrot and
  {M}ih\u ailescu)', \emph{J. Th\'eor. Nombres Bordeaux} \textbf{17} (2005),
  no.~1, 69--85.
\newblock \verb|http://jtnb.cedram.org/item?id=JTNB_2005__17_1_69_0|.

\bibitem{MR1484478}
W.~Bosma, J.~Cannon, and C.~Playoust,  `The {M}agma algebra system. {I}. {T}he
  user language', \emph{J. Symbolic Comput.} \textbf{24} (1997), no.~3-4,
  235--265.
\newblock \verb|http://dx.doi.org/10.1006/jsco.1996.0125|.

\bibitem{BT}
A.~Bremner and N.~Tzanakis, \emph{On the equation $y^2=x^6+k$}.
\newblock \verb|http://arxiv.org/abs/1001.3573|.

\bibitem{MR1144763}
J.~W.~S. Cassels, \emph{Lectures on elliptic curves}, in \emph{London
  Mathematical Society Student Texts} \textbf{24} (Cambridge University Press,
  Cambridge, 1991).

\bibitem{catalan}
E.~Catalan,  `Note extraite d'une lettre adress{\'e}e {\`a} l'editeur par {M}r.
  {E}. {C}atalan, {R}{\'e}p{\'e}titeur {\`a} l'{\'e}cole polytechnique de
  {P}aris', \emph{J. Reine Angew. Math.} \textbf{27} (1844), 192.
\newblock \verb|http://dx.doi.org/10.1515/crll.1844.27.192|.

\bibitem{MR0245500}
L.~E. Dickson, \emph{History of the theory of numbers. {V}ol. {II}:
  {D}iophantine analysis} (Chelsea Publishing Co., New York, 1966).

\bibitem{elkieshall}
N.~D. Elkies, \emph{List of integers {$x,y$} with {$x<10^{18}$}, {$0<|x^3-y^2|<
  x^{1/2}$}}.
\newblock \verb|http://www.math.harvard.edu/~elkies/hall.html|.

\bibitem{MR2548983}
G.~Everest and V.~Mah{\'e},  `A generalization of {S}iegel's theorem and
  {H}all's conjecture', \emph{Experiment. Math.} \textbf{18} (2009), no.~1,
  1--9.
\newblock \verb|http://projecteuclid.org/getRecord?id=euclid.em/1243430526|.

\bibitem{MR2045409}
G.~Everest, V.~Miller, and N.~Stephens,  `Primes generated by elliptic curves',
  \emph{Proc. Amer. Math. Soc.} \textbf{132} (2004), no.~4, 955--963
  (electronic).
\newblock \verb|http://dx.doi.org/10.1090/S0002-9939-03-07311-8|.

\bibitem{MR2041076}
G.~Everest, P.~Rogers, and T.~Ward,  `A higher-rank {M}ersenne problem', in
  \emph{Algorithmic number theory ({S}ydney, 2002)}, in \emph{Lecture Notes in
  Comput. Sci.} \textbf{2369}, pp.~95--107 (Springer, Berlin, 2002).
\newblock \verb|http://dx.doi.org/10.1007/3-540-45455-1_8|.

\bibitem{MR1464020}
J.~Gebel, A.~Peth{\"o}, and H.~G. Zimmer,  `Computing integral points on
  {M}ordell's elliptic curves', \emph{Collect. Math.} \textbf{48} (1997),
  no.~1-2, 115--136.
\newblock Journ{\'e}es Arithm{\'e}tiques (Barcelona, 1995).

\bibitem{MR1602064}
J.~Gebel, A.~Peth{\"o}, and H.~G. Zimmer,  `On {M}ordell's equation',
  \emph{Compositio Math.} \textbf{110} (1998), no.~3, 335--367.
\newblock \verb|http://dx.doi.org/10.1023/A:1000281602647|.

\bibitem{MR0323705}
M.~Hall, Jr.,  `The {D}iophantine equation {$x\sp{3}-y\sp{2}=k$}', in
  \emph{Computers in number theory ({P}roc. {S}ci. {R}es. {C}ouncil {A}tlas
  {S}ympos. {N}o. 2, {O}xford, 1969)}, pp.~173--198 (Academic Press, London,
  1971).

\bibitem{MR568909}
G.~H. Hardy and E.~M. Wright, \emph{An introduction to the theory of numbers}
  (The Clarendon Press Oxford University Press, New York, fifth ed., 1979).

\bibitem{MR2220098}
H.~A. Helfgott and A.~Venkatesh,  `Integral points on elliptic curves and
  3-torsion in class groups', \emph{J. Amer. Math. Soc.} \textbf{19} (2006),
  no.~3, 527--550 (electronic).
\newblock \verb|http://dx.doi.org/10.1090/S0894-0347-06-00515-7|.

\bibitem{MR1745599}
M.~Hindry and J.~H. Silverman, \emph{Diophantine geometry: An introduction}, in
  \emph{Graduate Texts in Mathematics} \textbf{201} (Springer-Verlag, New York,
  2000).
\newblock
  \verb|http://www.springer.com/mathematics/algebra/book/978-0-387-98975-4|.

\bibitem{MR2468477}
P.~Ingram,  `Multiples of integral points on elliptic curves', \emph{J. Number
  Theory} \textbf{129} (2009), no.~1, 182--208.
\newblock \verb|http://dx.doi.org/10.1016/j.jnt.2008.08.001|.

\bibitem{MR2076124}
P.~Mih{\u{a}}ilescu,  `Primary cyclotomic units and a proof of {C}atalan's
  conjecture', \emph{J. Reine Angew. Math.} \textbf{572} (2004), 167--195.
\newblock \verb|http://dx.doi.org/10.1515/crll.2004.048|.

\bibitem{jfm.44.0230.03}
L.~J. Mordell,  `The diophantine equation {$y^2-k=x^3$}', \emph{Proc. London
  Math. Soc.} \textbf{13} (1914), 60--80.
\newblock \verb|http://dx.doi.org/10.1112/plms/s2-13.1.60|.

\bibitem{mordell1914}
L.~J. Mordell,  `Indeterminate equations of the third and fourth degrees',
  \emph{Quart. J. of Pure and Applied Math.} \textbf{45} (1914), 170--186.

\bibitem{mordell1920}
L.~J. Mordell,  `A statement by {F}ermat', \emph{Proc. London Math. Soc.}
  \textbf{18} (1920), v.
\newblock \verb|http://dx.doi.org/10.1112/plms/s2-18.1.1-s|.

\bibitem{mordell1922}
L.~J. Mordell,  `On the rational solutions of the indeterminate equations of
  the third and fourth degrees', \emph{Proc. Cam. Phil. Soc.} \textbf{21}
  (1922), 179.

\bibitem{MR0186624}
D.~Mumford,  `A remark on {M}ordell's conjecture', \emph{Amer. J. Math.}
  \textbf{87} (1965), 1007--1016.
\newblock \verb|http://www.jstor.org/stable/2373258|.

\bibitem{MR992208}
J.~Oesterl{\'e},  `Nouvelles approches du ``th\'eor\`eme'' de {F}ermat',
  \emph{Ast\'erisque} (1988), no.~161-162, Exp.\ No.\ 694, 4, 165--186 (1989).
\newblock S{\'e}minaire Bourbaki, Vol. 1987/88.

\bibitem{MR0013386}
S.~S. Pillai,  `On the equation {$2\sp x-3\sp y=2\sp X+3\sp Y$}', \emph{Bull.
  Calcutta Math. Soc.} \textbf{37} (1945), 15--20.

\bibitem{reynolds}
J.~Reynolds, \emph{Extending {S}iegel's theory for elliptic curves} (Ph.D.
  thesis, Univ. of East Anglia, 2008).

\bibitem{siegel1926}
L.~Siegel,  `The integer solutions of the equation {$y^2=
  ax^n+bx^{n-1}+\cdots+k$}', \emph{J. Lond. Math. Soc.} \textbf{1} (1926),
  no.~2, 66--68.
\newblock \verb|http://dx.doi.org/10.1112/jlms/s1-1.2.66|.

\bibitem{siegel1929}
L.~Siegel,  `{\"U}ber einige {A}nwendungen diophantische {A}pproximationen', in
  \emph{Collected Works}, pp.~209--266 (Springer--Verlag, Berlin, 1966).

\bibitem{MR2514094}
J.~H. Silverman, \emph{The arithmetic of elliptic curves}, in \emph{Graduate
  Texts in Mathematics} \textbf{106} (Springer, Dordrecht, second ed., 2009).
\newblock \verb|http://dx.doi.org/10.1007/978-0-387-09494-6|.

\bibitem{MR0340175}
H.~M. Stark,  `Effective estimates of solutions of some {D}iophantine
  equations', \emph{Acta Arith.} \textbf{24} (1973), 251--259.
\newblock \verb|http://matwbn.icm.edu.pl/ksiazki/aa/aa24/aa2433.pdf|.
\newblock Collection of articles dedicated to Carl Ludwig Siegel on the
  occasion of his seventy-fifth birthday, III.

\bibitem{dioph}
P.~L. Tannery (ed.), \emph{Diophanti {A}lexandrini {O}pera omnia: cum {G}raecis
  commentariis} (Teubner, 1893--1895).
\newblock \verb|http://www.archive.org/stream/diophantialexan03plangoog|.

\bibitem{thue1909}
A.~Thue,  `{\"U}ber {A}nn{\"a}herungswerte {A}legbraischer {Z}ahlen', \emph{J.
  Reine Angew. Math.} \textbf{135} (1909), 284--305.
\newblock \verb|http://dx.doi.org/10.1515/crll.1909.135.284|.

\bibitem{MR0404136}
R.~Tijdeman,  `Some applications of {B}aker's sharpened bounds to {D}iophantine
  equations', in \emph{S\'eminaire {D}elange-{P}isot-{P}oitou (16e ann\'ee:
  1974/75), {T}h\'eorie des nombres, {F}asc. 2, {E}xp. {N}o. 24}, p.~7
  (Secr\'etariat Math\'ematique, Paris, 1975).
\newblock \verb|archive.numdam.org/article/SDPP_1974-1975__16_2_A3_0.pdf|.

\bibitem{pillai}
M.~Waldschmidt, \emph{Perfect powers: Pillai's works and their developments}.
\newblock \verb|http://arxiv.org/abs/0908.4031|.

\end{thebibliography}

\providecommand{\bysame}{\leavevmode\hbox to3em{\hrulefill}\thinspace}

\end{document}